\documentclass[letterpaper,12pt]{article}

\pdfoutput=1
\usepackage{graphicx}
\usepackage{chngpage}
\usepackage{booktabs}
\usepackage{amsmath}
\usepackage{multirow}
\usepackage{amsfonts}

\usepackage{hyperref}

\makeatletter
\renewcommand\section{\@startsection {section}{1}{\z@}%
{-3.5ex \@plus -1ex \@minus -0.2ex}%
{2.3ex \@plus 0.2ex}%
{\normalfont\normalsize\bfseries}}

\renewcommand\subsection{\@startsection{subsection}{2}{\z@}%
{-3.25ex \@plus -1ex \@minus -0.2ex}%
{1.5ex \@plus 0.2ex}%
{\normalfont\normalsize\bfseries}}

\def\@seccntformat#1{\csname the#1\endcsname.\quad}
\makeatother

\newcommand{\overbar}[1]{\mkern 4.3mu\overline{\mkern-4.3mu#1\mkern-1.5mu}\mkern 1.5mu}

\begin{document}

\setlength{\baselineskip}{4.5ex}

\noindent
\textbf{\Large Bayesian, classical and hybrid methods of inference}\\[1ex]
\textbf{\Large when one parameter value is special}\\[2ex]

\noindent
\textbf{Russell J. Bowater}\\
\emph{Independent researcher. Corresponding author. Contact via given email address or via the
website:
\href{https://sites.google.com/site/bowaterfospage}{sites.google.com/site/bowaterfospage}}
\\[3ex]
\textbf{Ludmila E. Guzm\'{a}n-Pantoja}\\
\emph{Assistant professor, Institute of Agro-industries, Technological University of the Mixteca,
Carretera a Acatlima Km.\ 2.5, Huajuapan de Le\'{o}n, Oaxaca, C.P.\ 69000, Mexico.}\\[2ex]

\noindent
\textbf{\small Abstract:}
{\small
This paper considers the problem of making statistical inferences about a parameter when a narrow
interval centred at a given value of the parameter is considered special, which is interpreted as
meaning that there is a substantial degree of prior belief that the true value of the parameter
lies in this interval.
A clear justification of the practical importance of this problem is provided.
The main difficulty with the standard Bayesian solution to this problem is discussed and, as a
result, a pseudo-Bayesian solution is put forward based on determining lower limits for the
posterior probability of the parameter lying in the special interval by means of a sensitivity
analysis.
Since it is not assumed that prior beliefs necessarily need to be expressed in terms of prior
probabilities, nor that post-data probabilities must be Bayesian posterior probabilities, hybrid
methods of inference are also proposed that are based on specific ways of measuring and
interpreting the classical concept of significance.
The various methods that are outlined are compared and contrasted at both a foundational level,
and from a practical viewpoint by applying them to real data from meta-analyses that appeared in
a well-known medical article.}

\pagebreak

\section{Introduction}

In statistical inference, examples abound where a particular value of a continuous population
parameter of interest merits special attention. A standard case could be when the value of the
parameter indicates that, with respect to the measurement scale concerned, there is no difference
between the performance of two treatments, e.g.\ the odds ratio or relative risk is equal to one,
or that there is no relationship between two variables, e.g.\ the correlation coefficient is
equal to zero.
However, while scientists may like to debate whether a continuous parameter is equal
or not to its `special' value, it is a question that is usually inadequately framed in
statistical terms, since, in general, it would be extremely unlikely in an applied context that the
parameter would be exactly equal to its special value.
Indeed, this is one of the reasons why hypothesis tests based on the sharp null hypothesis that a
parameter is equal to a given special value have received much criticism, see for example,
Anscombe~(1990), Edwards~(1965), Greenland~(2011), Oakes~(1990) and Royall~(1986). 

One exception to this argument could be in a study of paranormal activity, where the parameter of
interest equals a given value when such activity is absent, while any other possible value of the
parameter implies that such activity is present. Here, there is likely to be a substantial level of
prior belief that the parameter equals its special value. Therefore, it would seem appropriate that
an analysis of the data from such a study should take into account the special status of this
value, either explicitly by placing a large prior probability on this value in a Bayesian analysis,
or implicitly in the interpretation of the results of classical statistical methods.

While this type of scenario is rare, a slight relaxation of our main assumption allows us to
broaden our range of attention to include many more practical examples that are closely related
to this particular case. For this reason, let us define the special value of the parameter as being
special because a priori it is considered quite likely that the true value of the parameter lies
in a narrow interval centred at the special value. Clearly the original case is included in this
modified scenario if the width of the interval is allowed to be zero.

To give an example of this more general case, let us suppose that a test treatment, e.g.\ a new
type of antibiotic, is being compared with a control treatment, e.g.\ a standard type of
antibiotic, in a specific context, e.g.\ treating a certain type of infection.
Here, the parameter of interest will be assumed to be a measure of the effect of the treatment
relative to the control.
If, in many previous studies, treatments that are similar to the test treatment have shown effects
relative to the control that are very small or not distinguishable from zero, then there is likely
to be a substantial level of belief that the test treatment will also have an effect that is close
to zero relative to the control. The same would apply if there are no studies of similar
treatments, but many studies have already shown that the test treatment has an effect close to zero
relative to the control in various different contexts (e.g.\ treating different types of
infection) that are considered similar to the context that is currently under study (e.g.\ the
current interest is in treating a new type of infection.)

A standard way of estimating a treatment effect is to use a confidence interval. However, if there
was a substantial level of belief before carrying out a study that the true effect of the
treatment would lie in a narrow interval centred at the value of zero for the effect measure,
i.e.\ the interval $[-\varepsilon, \varepsilon]$, then the interpretation of a confidence interval
for the true effect will be complicated. In particular, if relatively little else was known a
priori about the true effect, then the probability that a confidence interval contains the true
effect will be felt to be smaller if it does not overlap with the interval $[-\varepsilon,
\varepsilon]$ in comparison to it fully containing this interval, despite both intervals having the
same confidence coefficient. Therefore, it is of interest to know whether non-standard methods can
be used to make more appropriate inferences about the treatment effect in this type of situation.

Addressing this issue and issues of a similar nature is the aim of the present paper.
More specifically, our main aim will be to examine ways in which inferences can be made about a
parameter $\theta$ when there is a substantial degree of prior belief that $\theta$ lies in a
narrow interval $[\theta_0-\varepsilon, \theta_0+\varepsilon]$ centred at a given special value
$\theta_0$, but little is known about the true value of $\theta$ when it is conditioned not to lie
in this special interval. Since it may be difficult to precisely specify how wide the special
interval for $\theta$ should be, it will be desirable that any inferences made about $\theta$ are
to some degree insensitive to the width of this interval, i.e.\ to the value of $2\varepsilon$.

It will not necessarily be assumed that prior belief about $\theta$ can be adequately represented
by a probability distribution over $\theta$, in this sense prior belief and prior probability will
be treated as distinct concepts. Also, although inferences about $\theta$ will be sought that
facilitate the construction of a probability distribution for $\theta$ that is valid after the data
has been observed, it will not be assumed that such a distribution must have been derived by the
usual Bayesian update of a prior to a posterior distribution. In this sense, an open view will be
taken as to what school of inference should be adopted to tackle the problem of interest.

In the following section, the standard Bayesian solution to this problem is reviewed, and the main
drawback of this solution is identified. To try to overcome this drawback, a two-step
Bayesian procedure is put forward in Section~\ref{sec2.2} that is based on performing a complete
sensitivity analysis over a very general class of prior densities for the parameter $\theta$. In
Section~\ref{sec3}, hybrid methods are proposed for tackling the general problem of interest that
are based on interpreting, in a specified way, measures of the classical concept of significance,
namely one-sided P values and what will be called Q values. The various methods of inference
presented in the paper are then applied to real data from a published set of meta-analyses in
Section~\ref{sec4}, and the final section of the paper contains a general discussion of the
advantages and disadvantages of these methods.
\vspace*{3ex}

\section{Bayesian approaches}
\vspace*{1ex}

\subsection{Standard Bayesian approach}
\label{sec2.1}

For the moment, let it be assumed that $\theta$ is the only unknown parameter in the sampling
model. Therefore, this model can be defined by the joint density of the data $x$ given $\theta$.
Let this sampling density be denoted as $f(x\,|\,\theta)$. To be able to use the Bayesian method to
make inferences about $\theta$, a prior density for $\theta$ needs to be specified. If this density
is defined to be $p(\theta)$, then it follows from Bayes' theorem that the density of $\theta$
conditioned on the data, i.e.\ the posterior density of $\theta$, is given by
\begin{equation}
\label{equ1}
p(\theta\,|\,x) = C f(x\,|\,\theta) p(\theta)
\end{equation}
where $C$ is a normalizing constant.

To take into account that there is a substantial degree of prior belief that $\theta$ lies in the
special interval $[\theta_0 - \varepsilon, \theta_0 + \varepsilon]$, a sizeable prior probability
$\alpha$ could be placed on the event of $\theta$ lying in this interval. Given that it has been
assumed that little would be known a priori about $\theta$ if it was conditioned not to lie in the
special interval, the prior distribution for $\theta$, under this condition, could be chosen to
have a large variance, and be a member of the family of all densities that are symmetric around
$\theta_0$ and monotonically decreasing as $|\theta-\theta_0|$ goes to infinity.
Let this family of distributions be denoted as $G_{\mbox{\scriptsize NIS}}$, where NIS is an
abbreviation for non-increasing and symmetric.

The complete prior density for $\theta$ is therefore defined as:
\begin{equation}
\label{equ2}
p(\theta) = \left\{
\begin{array}{ll}
\alpha h(\theta) & \ \ \mbox{if } \theta \in [\theta_0 - \varepsilon, \theta_0 +
\varepsilon],\vspace*{1ex}\\
(1-\alpha) g(\theta) & \ \ \mbox{otherwise}.
\end{array} \right.
\end{equation}
where $g(\theta) \in G_{\mbox{\scriptsize NIS}}$. Here $h(\theta)=1$ if $\varepsilon=0$, otherwise
$h(\theta)$ is any given prior density for $\theta$ conditional on $\theta \in
[\theta_0 - \varepsilon, \theta_0+\varepsilon]$. The prior density defined in equation~(\ref{equ2})
is updated to the posterior density for $\theta$ in the usual way by applying
equation~(\ref{equ1}), and a standard Bayesian solution to the problem of interest is thereby
obtained.

There is though a major drawback with this solution. The difficulty lies in the vagueness in
deciding how large the variance of the density $g(\theta)$ should be. For this conditional prior
density to be a good approximation to knowing very little about $\theta$ conditional on $\theta$
not lying in the interval $[\theta_0 - \varepsilon, \theta_0 + \varepsilon]$, we would like to
choose this variance to be as large as possible. Also, since any value chosen for this variance
will naturally be an imprecise value, we would like the resulting posterior distribution for
$\theta$ to be broadly insensitive to quite large changes in this variance.

However, it is well known that if, for any given sampling density $f(x\,|\,\theta)$ and any fixed
data set $x$, the variance of $g(\theta)$ is allowed to tend to infinity in way that ensures $\max
\{\hspace*{0.05em} g(\theta) : \theta \in \mathbb{R} \}$ tends to zero, then the posterior
probability that $\theta$ lies in the interval $[\theta_0 - \varepsilon, \theta_0 + \varepsilon]$
will tend to one (see Edwards, Lindman and Savage~1963 and Berger and Delampady~1987 for analysis
of related issues).
Therefore if, in these circumstances, the variance of $g(\theta)$ is chosen to be large enough, the
data $x$ would essentially be non-informative about whether or not $\theta$ lies in the interval
$[\theta_0 - \varepsilon, \theta_0 + \varepsilon]$.
This outcome is clearly unsatisfactory, and motivates the need to consider the kinds of alternative
approaches that will be discussed in the rest of this paper.
\vspace*{3ex}

\subsection{Two-step Bayesian approach}
\label{sec2.2}

Although, in many situations, a Bayesian sensitivity analysis may produce useful estimates of the
posterior probabilities of events of interest, it would appear to fail in the case just mentioned,
as the upper limit of the posterior probability of the event of $\theta$ lying in the interval
$[\theta_0-\varepsilon, \theta_0+\varepsilon]$, over all prior distributions under consideration,
would be equal to one for any given data set $x$.
Nevertheless, it may be worthwhile completing this Bayesian sensitivity analysis by determining a
lower limit for this posterior probability, as the value for this lower limit may prove to be an
acceptable rough estimate, or at least a useful underestimate, of the posterior probability for the
event in question.

Edwards, Lindman and Savage~(1963) carried out this kind of analysis under the assumption that
$\varepsilon=0$, implying that a prior probability of $\alpha$ is placed on the event of $\theta$
equalling $\theta_0$, and presented lower limits for the posterior probability of this event in
the case where the conditional prior density $g(\theta)$ belongs to the class of all possible
densities for $\theta$, and also in the case where it belongs to the class consisting only of
normal densities centred at the special value $\theta_0$.
Following on from this work, Berger and Sellke~(1987) derived these kind of lower limits,
under the assumption again that $\varepsilon=0$, but in the case where $g(\theta)$ is allowed to be
any density within the class of densities for $g(\theta)$ considered in the previous section, i.e.\
the class $G_{\mbox{\scriptsize NIS}}$.
In the two studies just mentioned, and in more detail in Berger and Delampady~(1987) and
Delampady~(1989), it is argued that allowing $\varepsilon$ to be positive but small, implying that
the special interval $[\theta_0-\varepsilon, \theta_0+\varepsilon]$ is indeed a narrow interval
rather than a single point, usually should not greatly affect the lower limits for the posterior
probability of $\theta$ lying in this interval calculated under the assumption that
$\varepsilon=0$.

The general conclusion of these studies is that these lower limits will be substantially above
zero, except if the height of the likelihood function for $\theta$ is very low at and around the
value $\theta_0$ relative to its maximum height.
Moreover, Berger and Delampady~(1987) made the case that, for any given data set $x$, these lower
limits will often be useful conservative approximations to the posterior probabilities for the
event of $\theta$ lying in $[\theta_0-\varepsilon, \theta_0+\varepsilon]$ that could result from
being able to determine a more accurate prior for $\theta$.

The authors of these studies, though, ignore the issue of whether at these lower limits, the form
of the prior density for $\theta$ will be acceptable in practice over all values of $\theta$.
For example, let us consider the case where $\varepsilon=0$ and the sampling density
$f(x\,|\,\theta)$ is defined by:
\[
f(x\,|\,\theta) = \prod^{n}_{i=1}\psi (x_i\,|\,\theta)
\]
where $x=(x_1, x_2, \ldots, x_n)$ and $\psi (x_i\,|\,\theta)$ is the normal density function with
unknown mean $\theta$ but known variance $\sigma^2$ evaluated at $x_i$.
In this case, at the lower limit for the posterior probability of the event that $\theta=\theta_0$
over all $g(\theta) \in G_{\mbox{\scriptsize NIS}}$, the density $g(\theta)$ will be a uniform
density conditioned not to lie in the interval $[\theta_0 - \varepsilon, \theta_0 + \varepsilon]$,
and it is possible that this density and, as a result, the prior density for $\theta$ will be zero
for all values of $\theta$ that are more than a short distance from $\theta_0$.
In fact, if the distance between the sample mean $\bar{x}$ and $\theta_0$ is less than 2.47 times
the standard error of the mean, i.e.\ $|\bar{x}-\theta_0| < 2.47\sigma/\sqrt{n}$, then the prior
density for $\theta$ will be zero for all values of $\theta$ that are further from $\theta_0$ than
the distance between the sample mean and $\theta_0$ plus at most just one standard error
$\sigma/\sqrt{n}$.
Given that the prior density for $\theta$ conditional on $\theta$ not lying in the interval
$[\theta_0 - \varepsilon, \theta_0 + \varepsilon]$ is supposed to represent a lack of prior
knowledge about $\theta$, this is clearly an unsatisfactory outcome.

A controversial, but perhaps acceptable, way of resolving this issue is to break down the inference
problem into two steps. In the first step, the posterior density for $\theta$ conditional on
$\theta$ not lying in the interval $[\theta_0-\varepsilon, \theta_0+\varepsilon]$ is determined
using a prior density for $\theta$ belonging to the class of densities 
$G_{\mbox{\scriptsize NIS}}$ that can be loosely described as being relatively flat or diffuse.
A critical aspect of this two-step procedure is that, by contrast to the behaviour of the full
posterior density for $\theta$ that was highlighted in Section~\ref{sec2.1}, this conditional
posterior density, which will be denoted as $p_{D}(\theta\,|\,\theta \notin [\theta_0-\varepsilon,
\theta_0+\varepsilon],x)$, is broadly insensitive to the choice made for the type of prior density
for $\theta$ under consideration.

To implement the second step of the method, we first need to define a class of densities for the
conditional prior density $g(\theta)$, which, for the purpose of giving an example, will again be
taken as being the class $G_{\mbox{\scriptsize NIS}}$, and also, an appropriate class of densities
$H$ for the prior density of $\theta$ conditional on $\theta$ lying in the interval $[\theta_0 - 
\varepsilon, \theta_0 + \varepsilon]$, i.e.\ a class of densities for $h(\theta)$. In this second
step, the posterior probability of the event of $\theta$ lying in $[\theta_0-\varepsilon,
\theta_0+\varepsilon]$ is then defined as being the minimum value of this posterior probability
over all $g(\theta) \in G_{\mbox{\scriptsize NIS}}$ and all $h(\theta) \in H$ with respect to a
given prior probability $\alpha$ assigned to this event.
Let this lower limit on the posterior probability in question be denoted as
$p_{L}(\theta \in [\theta_0-\varepsilon, \theta_0+\varepsilon]\,|\,x)$, and let the posterior
density for $\theta$ conditioned on $\theta$ lying in $[\theta_0-\varepsilon,
\theta_0+\varepsilon]$ that is obtained at the point where this lower limit is achieved be denoted
by $p_{L}(\theta\,|\,\theta \in [\theta_0-\varepsilon, \theta_0+\varepsilon],x)$.

We now have sufficient information to define the full posterior density for $\theta$ that
results from using this two-step Bayesian approach:
\[
p(\theta\,|\,x)\! =\! \left\{\!
\begin{array}{ll}
p_{L} (\theta\,|\,\theta \in [\theta_0-\varepsilon, \theta_0+\varepsilon],x) p_{L}(\theta \in
[\theta_0-\varepsilon, \theta_0+\varepsilon]\,|\,x) & \mbox{if }
\theta \in [\theta_0-\varepsilon, \theta_0+\varepsilon],\vspace*{1.5ex}\\
p_{D} (\theta\,|\,\theta \notin [\theta_0-\varepsilon, \theta_0+\varepsilon],x) (1 - p_{L}(\theta
\in [\theta_0-\varepsilon, \theta_0+\varepsilon]\,|\,x)) &
\mbox{otherwise}. \end{array} \right.
\]
This solution is not though a Bayesian solution in the true sense of the word as it is based on
the combination of two distinct choices for the prior density of $\theta$, i.e.\ the one used in
the first step of the method and the one used in the second step. Therefore, despite its intuitive
justification, the method is objectionable at a foundational level due to its incoherence from a
Bayesian perspective. Nevertheless, at least it may be considered as being a more viable way of
making inferences about $\theta$ than is offered by the standard Bayesian approach discussed in
Section~2.1.
\vspace*{3ex}

\section{Classical and hybrid approaches}
\label{sec3}
\vspace*{1ex}

\subsection{Interpretation of P values}

We will now consider ways of tackling the problem of interest that are based on the classical
concept of statistical significance, in particular, on the use of measures of significance such as
P values.
In what follows, $\hat{\theta}$ will be taken as being an estimator of $\theta$ that makes good
use of the information about $\theta$ that is contained in the sample, e.g.\ the maximum
likelihood estimator of $\theta$.
The estimator $\hat{\theta}$ will be assumed to be the test statistic that is used to
calculate P values with regard to hypotheses about $\theta$.
We will begin by establishing how P values will be interpreted in the methodology that will be
developed. This will be done by means of a simple example.

In this regard, let us imagine a scenario in which an individual who claims to have extrasensory
perception has correctly predicted the outcome of $y$ out of 20 independent Bernoulli trials in
succession, where the two outcomes of each trial both have a probability of 0.5 and $y$ is greater
than 10. It will be supposed that we wish to test the null hypothesis that the probability of a
correct prediction $p$ is equal to 0.5, i.e.\ that the individual does not have extrasensory
perception, against the alternative hypothesis that the probability $p$ is greater than 0.5, i.e.\
that he is using extrasensory perception to push his predictions in the right direction, or else,
he has cheated in some way. Under the assumption that the test statistic is $\hat{p}=y/20$, i.e.\
the maximum likelihood estimator of $p$, the P value for this problem, which clearly must be the
one-sided P value, is given by:
\[
\mbox{PV}(y) = \sum_{k=y}^{20} P(\mbox{$k$ correct predictions}) = \sum_{k=y}^{20} \binom{20}{k}
(0.5)^{20}
\]

The interpretation of this P value by any given analyst can be regarded as being dependent on
the analyst's opinion about the likely outcome of this experiment if it could be repeated under
exactly the same conditions. The repetition of the experiment will be treated as being only
hypothetical. Let the number of correct predictions made in this hypothetical future experiment be
denoted as $Y^*$. To be more specific, the key element in the interpretation of the P value will be
assumed to be the analyst's opinion regarding the plausibility of the following two
scenarios:\\[2ex]
Scenario 1: $P(Y^{*} \geq y)=\mbox{PV}(y)$, i.e.\ the probability that $Y^{*}$ is greater than or
equal to the observed number of correct predictions is equal to the P value.\\
Scenario 2: $P(Y^{*} \geq y)>\mbox{PV}(y)$, i.e.\ the probability that $Y^{*} \geq y$ is greater
than the P value.\\[2ex]
Clearly, Scenario~1 would be true if and only if the null hypothesis was true, while Scenario~2
would be true if and only if the alternative hypothesis was true. It could be argued therefore
that we are attempting to do nothing more than present an alternative way of defining the main
objective of a hypothesis test, i.e.\ to weigh up the plausibility of the null and alternative
hypotheses.

However, the interpretation of the P value via the comparison of the plausibility of these two
scenarios is useful, as it obliges the analyst to take into account both what he believed a priori
about the probability of a correct prediction, and the summary of information contained in the
data that is provided by the P value or, more precisely, the consistency that the data shows with
the null hypothesis as measured by the P value.
For example, if the analyst had a high degree of belief that the probability of a correct
prediction was 0.5 before the experiment was carried out, which of course would be a justifiable
belief given the context, then a value for PV(y) of 0.021, which corresponds to $y$ being equal to
15, may not be small enough to dissuade him from considering Scenario~1 as being much more
plausible than Scenario~2. However, if the same prior belief is combined with a value for PV(y) of
0.0002, which corresponds to $y$ being equal to 17, then the analyst may decide that the P value is
now small enough for him to regard Scenario~2 as being more plausible than Scenario~1.
In the rest of this paper, P values will be interpreted using this type of reasoning.
\vspace*{3ex}

\subsection{Hybrid method using one-sided P values}
\label{sec3.2}

Returning to the general problem of interest, let us assume that the observed value of the
estimator $\hat{\theta}$, as defined earlier, is such that if $\theta = \theta_0 - \varepsilon$
then the cumulative density function of $\hat{\theta}$, i.e.\
$F(\hat{\theta}\,|\,\theta_0 - \varepsilon)$, evaluated at this observed value, is less than or
equal to 0.5. As a result of this assumption, we will consider the case where the null and
alternative hypotheses are $H_0: \theta \geq \theta_0 - \varepsilon$ and
$H_1: \theta < \theta_0 - \varepsilon$ respectively.
For this case, $F(\hat{\theta}\,|\,\theta_0 - \varepsilon)$ is clearly the (one-sided) P value of
interest. Similar to what was discussed in the previous section, if $\hat{\theta}^{*}$ denotes the
value of $\hat{\theta}$ calculated on the basis of a future hypothetical sample of the same size as
the observed sample, then an analyst could try to interpret this P value by weighing up the
plausibility of the following two scenarios:\\[2ex]
Scenario 1: $P(\hat{\theta}^{*} \leq \hat{\theta}) \leq F(\hat{\theta}\,|\,\theta_0 -
\varepsilon)$, i.e.\ the probability that the hypothetical future estimate of $\theta$ is less than
or equal to the observed estimate of $\theta$ is less than or equal to the P value.\\
Scenario 2: $P(\hat{\theta}^{*} \leq \hat{\theta}) > F(\hat{\theta}\,|\,\theta_0 - \varepsilon)$,
i.e.\ the probability that $\hat{\theta}^{*} \leq \hat{\theta}$ is greater than the P value.\\[2ex]
Given the earlier definition of the estimator $\hat{\theta}$, it can easily be shown that, in
general, Scenario~1 would be true if and only if $\theta$ is greater than or equal to
$\theta_0 - \varepsilon$, while Scenario~2 would be true if and only if $\theta$ is less than
$\theta_0 - \varepsilon$.

To weigh up the plausibility of these two scenarios, the analyst is obliged to take into account
both what he felt about the likeliness of the events $\{ \theta \geq \theta_0 - \varepsilon \}$ and
$\{ \theta < \theta_0 - \varepsilon \}$ before the data was observed, and the size of the P value.
In particular, since to be in the context of interest, it needs to be assumed that the analyst had
a substantial degree of prior belief that $\theta$ lay in the special interval $[\theta_0-
\varepsilon, \theta_0+\varepsilon]$, but also would have known little about $\theta$ if it was
conditioned not to lie in this interval, it would be expected that quite a small P value would be
required in order to dissuade him from believing that Scenario 1 is more plausible than Scenario~2.

If the analyst is prepared to express his opinion about the plausibility of Scenarios~1 and~2 by
assigning probabilities to these two scenarios, then he will of course have effectively determined
probabilities for the events $\{ \theta \geq \theta_0 - \varepsilon \}$ and $\{ \theta < \theta_0 -
\varepsilon \}$ that are applicable after the data has been observed.
As mentioned in the Introduction, it will not be assumed that the validity of these post-data
probabilities depends on how acceptable they would be if they were treated as posterior 
probabilities derived using the usual Bayesian updating rule.
However, the procedure based on using one-sided P values that has been put forward does not
naturally lead us to a point where we would be able to determine a post-data density for $\theta$
over the whole of the real line, or even, a post-data probability for the event of $\theta$ lying
in the special interval $[\theta_0-\varepsilon, \theta_0+\varepsilon]$.

For this reason, we will specify two types of conditional prior density for $\theta$ that are
generally consistent with the post-data probabilities for the events $\{ \theta \geq \theta_0 -
\varepsilon \}$ and $\{ \theta < \theta_0 - \varepsilon \}$ that are produced by the above
procedure.
In particular, with the same justification as given in Section~\ref{sec2.1}, let us again specify
the prior density for $\theta$ conditional on $\theta$ not lying in the interval
$[\theta_0-\varepsilon, \theta_0+\varepsilon]$ as being a relatively diffuse density from the
class of densities $G_{\mbox{\scriptsize NIS}}$.

Without the need for further assumptions, we are now able to deduce what would be the post-data
probability of the event of $\theta$ lying in $[\theta_0-\varepsilon, \theta_0+\varepsilon]$. This
probability is
\begin{equation}
\label{equ3}
P(\theta \in [\theta_0-\varepsilon, \theta_0+\varepsilon]\,|\,x) = \gamma - \lambda (1-\gamma)
\end{equation}
where $\gamma$ is the post-data probability assigned to the event $\{ \theta \geq \theta_0 -
\varepsilon \}$ using the procedure based on the one-sided P value outlined above, and where,
conditional on $\theta$ not lying in $[\theta_0-\varepsilon, \theta_0+\varepsilon]$, the value
$\lambda$ is the posterior probability of the event $\{ \theta > \theta_0 + \varepsilon \}$ divided
by the posterior probability of the event $\{\theta < \theta_0 - \varepsilon \}$, which, using the
same notation as in Section~\ref{sec2.2}, can be expressed as
\[
\lambda = \frac{p_{D}(\theta > \theta_0 + \varepsilon \,|\, \theta \notin [\theta_0 -
\varepsilon, \theta_0 + \varepsilon], x)}{p_{D}(\theta < \theta_0 - \varepsilon \,|\, \theta
\notin [\theta_0 - \varepsilon, \theta_0 + \varepsilon], x)}
\]
Also if, as was done in Section~\ref{sec2.1}, we specify the prior density for $\theta$ conditional
on $\theta$ lying in the interval $[\theta_0-\varepsilon, \theta_0+\varepsilon]$ as being any given
probability density $h(\theta)$ that is restricted to lie on this interval, then it can be deduced
that the post-data density for $\theta$ over the whole of the real line would be
\vspace*{1ex}
\begin{equation}
\label{equ4}
P(\theta\,|\,x) = \left\{
\begin{array}{ll}
(1-\gamma)p_{D}(\theta\,|\,\{ \theta < \theta_0 - \varepsilon \}, x) &
\mbox{if } \theta < \theta_0 - \varepsilon,\vspace*{1.5ex}\\
(\gamma-\lambda(1-\gamma))P(\theta\,|\, \theta \in [\theta_0 - \varepsilon, \theta_0
+ \varepsilon], x)
& \mbox{if } \theta \in [ \theta_0 - \varepsilon, \theta_0 + \varepsilon ],\vspace*{1.5ex}\\
\lambda(1-\gamma)p_{D}(\theta\,|\,\{ \theta > \theta_0 + \varepsilon \}, x) &
\mbox{if } \theta > \theta_0 - \varepsilon\hspace*{20em}
\end{array}
\right.
\end{equation}
where all conditional posterior densities on the right-hand side of this expression are
directly derived from their corresponding conditional prior densities.

Similar to one of the characteristics of the two-step Bayesian method, and in contrast to the
posterior density for $\theta$ that was highlighted in Section~\ref{sec2.1}, the post-data density
$P(\theta\,|\,x)$ defined in equation~(\ref{equ4}) is largely insensitive to the choice made for
the conditional prior density $P(\theta\,|\, \theta \notin [ \theta_0 - \varepsilon, \theta_0 +
\varepsilon ])$ when it has already been decided that the form of this prior density must be as
specified above. Nevertheless, this conditional prior density will not be consistent with the
post-data probability $\gamma$ assigned to the event of $\theta$ lying in the interval
$[\theta_0-\varepsilon, \theta_0+\varepsilon]$ unless $\gamma$ satisfies the following condition:
\begin{equation}
\label{equ6}
\gamma \geq \frac{\lambda}{1+\lambda}
\end{equation}
If $\gamma$ does not satisfy this condition then the post-data probability
$P(\theta \in [\theta_0-\varepsilon, \theta_0+\varepsilon]\,|\,x)$ defined in equation~(\ref{equ3})
will be negative. In practice, this condition though may provide the analyst with a useful
reference point rather than an obstacle in choosing an appropriate value for $\gamma$.

Although, at the start of this section, the assumption was made that the cumulative probability
$F(\hat{\theta}\,|\,\theta_0 - \varepsilon)$ is less than or equal to 0.5, it should be clear how
a similar method could be applied to derive a post-data density for $\theta$ in the case where this
cumulative probability calculated under the assumption that $\theta = \theta_0 + \varepsilon$,
i.e.\ $F(\hat{\theta}\,|\,\theta_0 + \varepsilon)$, is greater than or equal to 0.5.
In particular, it should be clear that it would be sensible to take the one-sided P value of
interest as being $1 - F(\hat{\theta}\,|\,\theta_0 + \varepsilon)$ in this case, which corresponds
to the null and alternative hypotheses being $H_0: \theta \leq \theta_0 + \varepsilon\,$ and $H_1:
\theta > \theta_0 + \varepsilon$ respectively.

What to do, on the other hand, in the rare case where $F(\hat{\theta}\,|\,\theta_0 - \varepsilon)$
is greater than 0.5 and $F(\hat{\theta}\,|\,\theta_0 + \varepsilon)$ is less than 0.5 is more of a
debating point. In these circumstances the estimator $\hat{\theta}$, though, would really not
provide any grounds for having less belief that the parameter $\theta$ lies in the special interval
$[\theta_0 - \varepsilon, \theta_0 + \varepsilon]$ than before the data was observed. Therefore, in
practice, an adequate solution to this problem may be to assign the same post-data probability to
the event $\theta \in [\theta_0 - \varepsilon, \theta_0 + \varepsilon]$ as the prior probability
that would have been given to this event.

Of course, it could be argued that since, in this case, the estimator $\hat{\theta}$ appears to be
corroborating the hypothesis that $\theta$ lies in the interval $[\theta_0 - \varepsilon, \theta_0
+ \varepsilon]$, the post-data probability that is assigned to this hypothesis should be greater
than the prior probability given to this hypothesis. However, the degree of corroboration is likely
to be very small in practice, and it could also be very difficult to precisely specify exactly by
how much the former probability should be greater than the latter.
\vspace*{3ex}

\subsection{Hybrid method using Q values}
\label{sec3.3}

Given that two-sided P values are usually the significance measure of choice when the null
hypothesis asserts that the parameter under consideration is equal to a given special value, and
the alternative hypothesis is the opposite of this hypothesis, we will now try to develop a method
similar to the one proposed in the previous section, but based on a measure of significance that
results from extending the concept of a two-sided P value to the general problem of interest.
For want of a better name, this measure of significance will be called a Q value, which should not
be confused with any other meanings given to the same term. To illustrate the concept of a Q value
and, more specifically, to show how such a measure can be used in the context being analysed, let
us for the moment consider the case where the parameter $\theta$ is the population mean $\mu$ with
its special interval denoted by $[\mu_0 - \varepsilon, \mu_0 + \varepsilon]$, and where the
sampling distribution of the data $x$ given $\mu$ is such that the data values independently follow
the same normal distribution with known variance $\sigma^2$.

Under the assumption that the test statistic is the sample mean $\bar{x}$ of $n$ observations, the
Q value that corresponds to the null hypothesis that $\mu=\mu_*$ will be defined as:
\vspace*{0.5ex}
\begin{equation}
\label{equ5}
q(\mu_*) = \left\{
\begin{array}{ll}
1 - F(\bar{x}\,|\,\mu_*) + F(2\mu_0 - \bar{x}\,|\,\mu_*) &
\ \ \mbox{if } \bar{x} \geq \mu_0,\vspace*{1.5ex}\\
1 - F(2\mu_0 - \bar{x}\,|\,\mu_*) + F(\bar{x}\,|\,\mu_*) &
\ \ \mbox{if } \bar{x} < \mu_0
\end{array}
\right.
\vspace*{0.5ex}
\end{equation}
where $F(w\,|\,\mu_*)$ is the cumulative normal density with mean $\mu_*$ and variance $\sigma^2/n$
evaluated at $w$.
Therefore, this Q value can be described as being the sum of the area under the normal density with
mean $\mu_*$ and variance $\sigma^2/n$ that lies to the left of $\mu_0 - d$ and to the right of
$\mu_0 + d$, where the distance $d$ is equal to $|\bar{x}-\mu_0|$. It can be seen that in the
case where $\mu_* = \mu_0$, the Q value is equivalent to the usual two-sided P value for the null
hypothesis concerned, i.e.\ the null hypothesis that $\mu = \mu_0$. Figure~1 illustrates how a Q
value is calculated in the case where $\mu_0=0$, in particular the Q value is equal to the total
shaded area under the curve in this figure.

\begin{figure}[tb]
\begin{flushleft}
{\small \bf Figure 1. Calculation of Q value in the case where $\mu_0=0$.}\\[0.5ex]
{\footnotesize Q value = total shaded area. The interval $(-\varepsilon, \varepsilon)$ has only
been included for illustrative purposes.}
\vspace*{-2ex}
\end{flushleft}
\begin{center}
\includegraphics[width=5.75in]{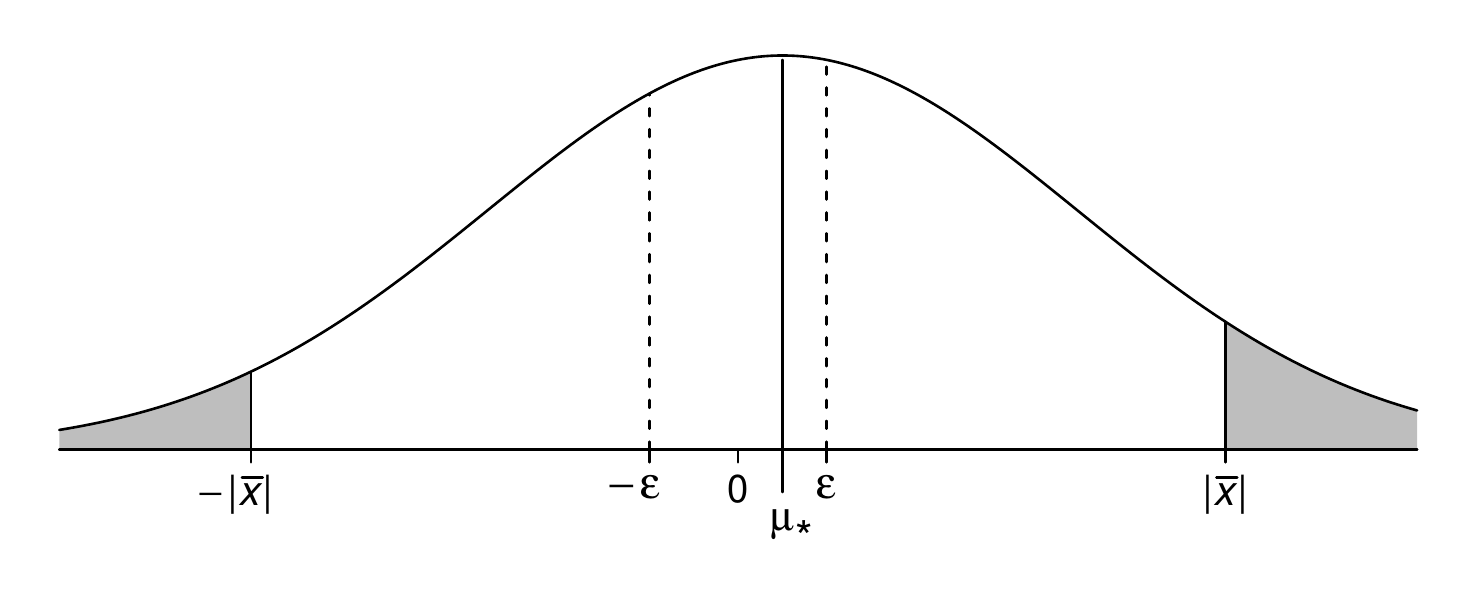}
\end{center}
\vspace*{-2ex}
\end{figure}

Due to the symmetry of the normal density, the Q value as defined in equation~(\ref{equ5}) has the
following property:
\[
q(\mu_0 + \varepsilon) = q(\mu_0 - \varepsilon)
\]
i.e.\ the Q value is equal to the same value at the two limits of the special interval
$[\mu_0-\varepsilon, \mu_0 + \varepsilon]$, and also the property:
\[
\begin{array}{ll}
q(\mu') \leq q(\mu_0 + \varepsilon) & \ \ \forall\, \mu' \in [\mu_0-\varepsilon, \mu_0 +
\varepsilon]\vspace*{1ex}\\
q(\mu') > q(\mu_0 + \varepsilon) & \ \ \forall\, \mu' \notin [\mu_0-\varepsilon, \mu_0 +
\varepsilon]
\end{array}
\]
For this reason, if $\overbar{X}^{*}$ denotes the mean of a future hypothetical sample of $n$
data values, then an analyst could try to make inferences about whether or not the
parameter $\mu$ lies in the interval $[\mu_0-\varepsilon, \mu_0 + \varepsilon]$ by weighing up the
plausibility of the following two scenarios:\\[2ex]
Scenario 1: $P(|\overbar{X}^{*}-\mu_0| > |\bar{x}-\mu_0|) \leq q(\mu_0 + \varepsilon)$, i.e.\ the
probability that $\overbar{X}^{*}$ is further away from $\mu_0$ than the distance between the
observed sample mean and $\mu_0$ is less than or equal to the Q value when $\mu_*=\mu_0 +
\varepsilon$.\\
Scenario 2: $P(|\overbar{X}^{*}-\mu_0| > |\bar{x}-\mu_0|) > q(\mu_0 + \varepsilon)$, i.e.\ the
opposite of the relationship that defines Scenario~1.\\[2ex]
It can easily be shown that Scenario 1 would be true if and only if $\mu$ lies in the special
interval $[\mu_0-\varepsilon, \mu_0 + \varepsilon]$, while Scenario 2 would be true if and only if
$\mu$ does not lie in this interval.

To weigh up the plausibility of these two scenarios, the analyst is obliged to take into account
both what he felt about the likeliness of $\mu$ lying in the interval $[\mu_0-\varepsilon,
\mu_0 + \varepsilon]$ before the data was observed, and the size of the Q value when
$\mu_*=\mu_0 + \varepsilon$.
The effect of a small value for $q(\mu_0 + \varepsilon)$ in cases where $\varepsilon>0$ should be
similar to the effect of a small value for the two-sided P value when the null hypothesis is
$\mu=\mu_0$, that is, it should generally lead to a lowering in the degree of belief that $\mu$
lies in the interval $[\mu_0 - \varepsilon, \mu_0 + \varepsilon]$ in comparison to what was felt
about the likeliness of $\mu$ lying in this interval before the data observed.
Indeed, since it has been assumed that $\varepsilon$ must be small, the value
$q(\mu_0 + \varepsilon)$ should generally be well approximated by the two-sided P value for the
null hypothesis that $\mu=\mu_0$, to the extent that if $\varepsilon=0$, then this
two-sided P value will be of course equal to the Q value in question.

If the analyst expresses what he believes about the plausibility of Scenarios~1 and~2 by assigning
probabilities to these two scenarios, then he will have effectively determined post-data
probabilities for the event of $\mu$ lying in the special interval $[\mu_0-\varepsilon, \mu_0 +
\varepsilon]$ and the complement of this event.
Assuming that this has been done, let us suppose that the prior density for $\mu$ conditional on
$\mu$ lying in the interval $[\mu_0-\varepsilon, \mu_0 + \varepsilon]$ and the prior density for
$\mu$ conditional on this not being the case are defined in the same way as in
Section~\ref{sec3.2}, with of course, the mean $\mu$ taking the place of the general parameter
$\theta$.
It is clear that any choice for these two conditional prior densities will be consistent with the
assignment of any post-data probability to the event of $\mu$ lying in $[\mu_0-\varepsilon,
\mu_0 + \varepsilon]$, in the sense that the integral of the resulting post-data density for
$\mu$ over all values of $\mu$ will be guaranteed to equal one.
If $\beta$ denotes the value assigned to the post-data probability of $\mu$ lying in
$[\mu_0-\varepsilon, \mu_0 + \varepsilon]$, then it can be deduced from the assumptions already
made that the post-data density for $\mu$ over all values of $\mu$ is given by
\vspace*{1ex}
\[
P(\mu\,|\,x) = \left\{
\begin{array}{ll}
\beta P(\mu\,|\,\mu \in [\mu_0 - \varepsilon, \mu_0 + \varepsilon], x) &\ \
\mbox{if } \mu \in [\mu_0 - \varepsilon, \mu_0 + \varepsilon],\vspace*{1.5ex}\\
(1 - \beta)p_{D}(\mu\,|\, \mu \notin [\mu_0 - \varepsilon, \mu_0 + \varepsilon], x) &\ \
\mbox{if } \mu \notin [\mu_0 - \varepsilon, \mu_0 + \varepsilon]
\end{array}
\right.
\vspace*{1ex}
\]
where the two conditional posterior densities on the right-hand side of this expression are
directly derived from their corresponding conditional prior densities.
As was the case for the methods proposed in Sections~\ref{sec2.2} and~\ref{sec3.2}, the post-data
density $P(\theta\,|\,x)$, or in the present context $P(\mu\,|\,x)$, is largely insensitive to the
choice made for the conditional prior density $P(\theta\,|\, \theta \notin [ \theta_0 -
\varepsilon, \theta_0 + \varepsilon ])$ when it has already been decided that this prior density
must take the form that has been specified.

Although up to now, we have only considered a case in which the sampling distribution of the
test statistic is normal, and in which the variance of this distribution is constant over different
values of the parameter of interest, the Q value method that has been outlined can be easily
extended to deal with cases in which this assumption of normality and/or this assumption of a
constant variance are inappropriate.
Nevertheless, this method is not as widely applicable as the one-sided P value method outlined in
the previous section.
In particular, it will generally be difficult to apply this method in a satisfactory manner when
the sampling distribution of the test statistic is substantially skewed. 

The advantage of the proposed method is that it facilitates the direct assignment of a post-data
probability to the event of the parameter $\theta$ lying in its special interval $[\theta_0 -
\varepsilon, \theta_0 + \varepsilon]$, instead of achieving this indirectly via first assigning a
post-data probability to the event that $\theta$ lies in $[\theta_0 - \varepsilon, \infty)$ or to
the event that $\theta$ lies in $(-\infty, \theta_0 + \varepsilon]$, as is the case for the
one-sided P value method. Also, it is possible to apply the method adequately in an approximate
sense in situations where it can not be applied exactly, as will be demonstrated in the next
section.
\vspace*{3ex}

\section{Methods applied to real data}
\label{sec4}
\vspace*{1ex}

\subsection{Overview}

We will now apply the methods that were outlined in previous sections to real data.
In particular, we will consider the set of meta-analyses presented in a famous paper by Andraws,
Berger and Brown that appeared in the Journal of the American Medical Association in 2005 (Andraws,
Berger and Brown 2005). These meta-analyses combine the results of randomized controlled trials
that studied the effect of antibiotic treatment as an addition to standard medical therapy in
patients diagnosed with coronary artery disease.

The main focus of our attention will be on the meta-analysis that relates to the effect of this
treatment on the occurrence of heart attacks and unstable angina, i.e.\ acute coronary syndromes.
This meta-analysis is presented in Figure~4 of the paper in question. As was done in this earlier
work, we will measure the effect of the treatment relative to the control using the odds ratio.
Also, as is standard practice, we will assume that the sampling distribution of the logarithm of
the odds ratio can be satisfactorily approximated by a normal distribution, and that the usual
estimate of the variance of this distribution can be generally treated as though it is the true
value of this variance, see Woolf~(1955) for more details.

The results of the four largest studies in the aforementioned meta-analysis are reproduced in the
first four rows of Figure~2 of the present paper. Combining the odds ratios of these four studies
using the random effects method of DerSimonian and Laird~(1986) produces a 95\% confidence interval
for the overall odds ratio of $(0.890, 1.100)$, as indicated in the 5th row of Figure~2, which on
the log scale is an interval of $(-0.116, 0.096)$.
Therefore, we would feel entitled to have a high degree of belief that the logarithm of the true
odds ratio for these studies would lie in the interval $[-0.1,0.1]$, which is clearly quite a
narrow interval centred at the null effect. Given that it will be assumed in what follows that the
four studies being considered represent all the information available `a priori' about the medical
question of interest, this interval will be treated as the special interval for the log odds ratio.
Since the log odds ratio will be designated as the parameter of interest $\theta$, it is of course
required that $\theta_0=0$ and $\varepsilon=0.1$ when this special interval is expressed using
earlier notation.

\begin{figure}[tb]
\vspace*{0.5ex}
\begin{flushleft}
{\small \bf Figure 2. Forest plot of the studies used in the data analysis.}\\[0.5ex]
{\footnotesize ACS = Acute coronary syndromes, CI = Confidence interval, Combined = the preceding
studies combined, MI = Myocardial infarction, OR Est. = Estimate of the odds ratio}
\vspace*{-1ex}
\end{flushleft}
\begin{adjustwidth}{-1in}{-1in}
\centering
\includegraphics[width=7in]{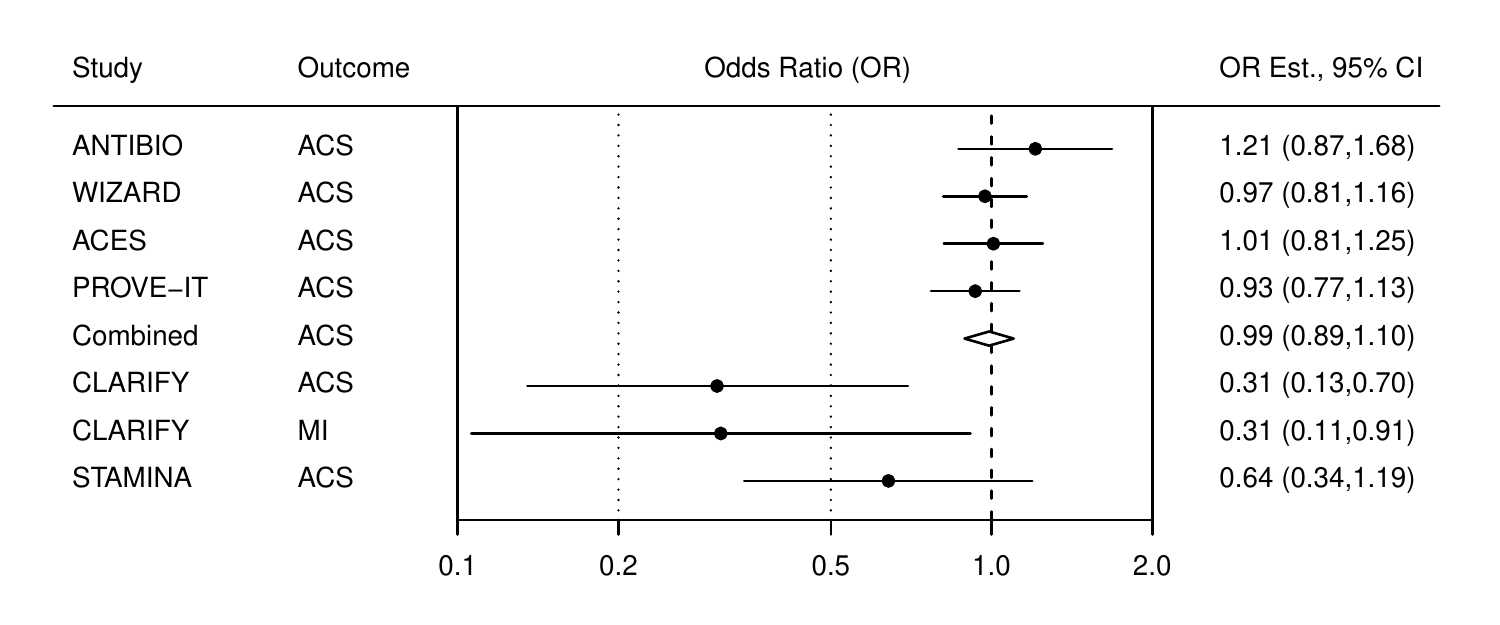}
\end{adjustwidth}
\vspace*{1ex}
\end{figure}

We will now attempt to make inferences about the log odds ratio $\theta$ on the basis of the
results of various possible new studies of the phenomenon in question.
\vspace*{3ex}

\subsection{First new study}
\label{sec4.2}

To begin with, let us imagine that the new data that needs to be analysed are the event rates with
respect to the same medical outcomes for the CLARIFY study, which is a relatively small clinical
trial that is also included in the meta-analysis that has just been highlighted. This study is
similar, but notably distinct, from the four larger studies mentioned above. For example, the
antibiotic used in this study is clarithromycin, which is not one of the types of antibiotic used
in these larger studies. Such differences may adequately justify why it is necessary to carry out
such a new study, but the similarities with the other studies are likely to imply the presence of a
substantial level of prior belief that the true value of the log odds ratio (log OR) for this study
will lie in the interval $[-0.1,0.1]$.

The event counts for this study can be found in the legend of Table~1 of the present paper.
The sample odds ratio for these event counts is 0.306 and, under the assumptions that were
outlined above, the 95\% confidence interval for the odds ratio in this study is $(0.135, 0.697)$,
as indicated in row~6 of Figure~2.
Clearly this interval does not contain the null effect, i.e.\ odds ratio $= 1$, and is indeed a
substantial distance away from the null effect. Under the same assumptions, the 95\% central
posterior interval (or credible interval) for the odds ratio has limits that are equal to the
limits of this interval if a flat improper prior is placed over log OR. Information relating to the
use of this prior distribution is given in the first row of Table~1.

\begin{table}[tb]
\begin{flushleft}
{\small {\bf Table 1. Analysis of data from the CLARIFY study for acute coronary
syndromes.}}\\[0.5ex]
{\footnotesize Data summary: $n_t$ (no.\ of patients in the treatment group) \hspace{-0.1em}= 
\hspace{-0.1em}74, $n_c$ (no.\ of patients in the control group) \hspace{-0.1em}=
\hspace{-0.1em}74, $e_t$ (no.\ of events in the treatment group) \hspace{-0.1em}=
\hspace{-0.1em}10, $e_c$ (no.\ of events in the control group) \hspace{-0.1em}=
\hspace{-0.1em}25.}\\
{\footnotesize Key: OR = odds ratio.}
\end{flushleft}
\renewcommand{\arraystretch}{1.1}
\begin{adjustwidth}{-1in}{-1in}
\begin{center}
{\footnotesize
\begin{tabular}{|l|l|l|l|l|l|}
\hline
Method & Prior probability & Significance & Post-data & Post-data & 95\% central\\
& of log OR $\in$ & measure & probability of & probability of & post-data\\
& $[-0.1,0.1]$ & & log \hspace*{-0.1em}OR $\!\in\! [-0.1,0.1]$ & log OR $\geq -0.1$
& interval for OR\\
\hline
Bayesian with flat & \multirow{2}{*}{N/a} & \multirow{2}{*}{N/a} & \multirow{2}{*}{0.004} &
\multirow{2}{*}{0.005} & \multirow{2}{*}{(0.135, 0.697)}\\
improper prior & & & & & \\
\hline
\multirow{2}{*}{Two-step Bayesian}
& 0.5* & \multirow{2}{*}{N/a} & 0.064 & 0.065 & (0.136, 0.992)\\
& 0.8* & & 0.214 & 0.215 & (0.141, 1.060)\\
\hline
\multirow{3}{28mm}{One-sided P value hybrid} & \multirow{3}{*}{N/a} &
\multirow{3}{23mm}{One-sided P value = 0.0049} & 0.049 & 0.05* & (0.136, 0.971)\\
& & & 0.019 & 0.02* & (0.135, 0.813)\\
& & & 0.099 & 0.01* & (0.135, 0.725)\\
\hline
\multirow{2}{28mm}{Q value hybrid} & \multirow{2}{*}{N/a} & \multirow{2}{23mm}{Q value = 0.0060}
& 0.05* & 0.051 & (0.136, 0.973)\\
& & & 0.01* & 0.011 & (0.135, 0.732)\\
\hline
\end{tabular}}
\end{center}
\end{adjustwidth}
\vspace*{-2ex}
\begin{flushleft}
{\footnotesize The asterisks (*) indicate probabilities determined by the analyst.}
\end{flushleft}
\vspace*{2ex}
\end{table}

Let us now analyse the same data set using the two-step Bayesian approach outlined in
Section~\ref{sec2.2}.
We will apply this method as it was specified in this earlier section with the prior density for
$\theta$ conditional on $\theta$ lying in the special interval $[-0.1,0.1]$, i.e.\ the density
$h(\theta)$, assumed to be a uniform density over this interval implying that the class of
densities $H$ contains just this single density.
Making this assumption would seem to be reasonable, and would not be expected to have major
practical consequences given the narrowness of the interval concerned. 
It will also be assumed that, in the first step of the two-step Bayesian approach, the prior
density for log OR conditional on log OR not lying in the interval $[-0.1,0.1]$, which,
according to the way this method was defined, is required to be a relatively diffuse density
from the class of densities $G_{\mbox{\scriptsize NIS}}$, is adequately approximated by a flat
improper prior over all values for log OR not lying in this interval.

Given that it has been argued that a substantial level of prior belief should be placed on the
event of log OR lying in the interval $[-0.1,0.1]$, it would seem appropriate that a prior
probability of say 0.5 or 0.8 is assigned to this event. Information relating to the application
of the two-step Bayesian method to the data set under consideration, using these two prior
probabilities, can be found in rows~2 and~3 of Table~1. It can be seen therefore that, for these
prior probabilities, the 95\% central post-data (or posterior) intervals for the odds ratio are
very different from the corresponding 95\% confidence interval given in row~1 of the table.
In particular, they substantially overlap with the interval $[0.905, 1.105]$, which corresponds to
the special interval of $[-0.1,0.1]$ for log OR. Nevertheless, the posterior probabilities of
log OR lying in the interval $[-0.1, 0.1]$ are substantially less than the prior probabilities
assigned to this event, as is indicated in column~4 of the table.

If the null and alternative hypotheses are $H_0$: log OR $\geq -0.1$ and $H_1$: log OR $< -0.1$
respectively, then the one-sided P value is 0.0049. Therefore, the two scenarios that need to be
weighed up in order to implement the P value method described in Section~\ref{sec3.2} are as
follows:\\[1.5ex]
Scenario 1: $P(\hspace*{0.2em}\mbox{sample OR*} < \mbox{observed odds ratio} =
0.306\hspace*{0.2em}) \leq 0.0049$\\
Scenario 2: $P(\hspace*{0.2em}\mbox{sample OR*} < 0.306\hspace*{0.2em}) > 0.0049$\\[1.5ex]
where sample OR* is the sample odds ratio in a future hypothetical sample with the number of
patients in the treatment and control groups being the same as in the observed sample.
Under the assumptions that have already been made about the distribution of the sample value of
log OR, Scenario~1 would be true if and only if log OR $\geq -0.1$, while Scenario~2 would be
true if and only if log OR $< -0.1$.

For the current data, the condition in equation~(\ref{equ6}) implies that the minimum probability
that should be assigned to Scenario 1 is 0.001. Furthermore, it can be seen from row~1 of Table~1
that if a flat improper prior is placed over all values of log OR, then a posterior probability of
0.005 should be assigned to the event $\{\hspace*{0.2em}\mbox{log OR} \geq -0.1\hspace*{0.05em}\}$.
Given that there was in fact a substantial degree of prior belief that log OR lay in $[-0.1,0.1]$,
we would expect therefore that the probability assigned to Scenario~1 would be substantially higher
than 0.005.

Clearly the assignment of probabilities to the two scenarios presented above depends on the
subjective judgement of the person analysing the data. In rows~4 to~6 of Table~1, various values
are given for the probability that could be assigned to Scenario~1 or equivalently, under the
assumptions that have been made, the probability that log OR $\geq -0.1$. The three values proposed
have been chosen to represent a high value (0.05), a medium value (0.02) and a low value
(0.01) for this probability. Nevertheless, all these values have been chosen to represent a
relatively small probability, since the small one-sided P value of 0.0049 would seem to heavily
disfavour the possibility of log OR lying in the interval $[-0.1,0.1]$.

The 95\% central post-data intervals for the odds ratio in rows~4 to~6 of Table~1 have been
calculated under the assumption that the prior density for log OR conditional on log OR lying in
the interval $[-0.1,0.1]$ is uniform over this interval and, conditional on log OR not lying in
this interval, takes the form of a flat improper density over all possible values of log OR.
It can be seen that these post-data intervals for the odds ratio all have upper limits
that are lower than the upper limits of the corresponding intervals for the two-step Bayesian
approach given in rows~2 and~3 of this table. However, as is the case with the two-step
Bayesian approach, the upper limits of these intervals are clearly distinct from the upper limit of
the 95\% confidence interval for the odds ratio given in row~1 of the table.

In Section~\ref{sec3.3}, it was described essentially, how, with respect to a specific null
hypothesis about a parameter $\theta$, the Q value for any given estimator of $\theta$ that has a
normal distribution can be determined. Placing this definition in the context of the present type
of data analysis implies that the Q value for the data set under consideration that corresponds to
the null hypothesis that log OR $=-0.1$ is 0.006. 
Therefore, the two scenarios that need to be weighed up in order to implement the Q value method
outlined earlier are as follows:\\[1.5ex]
Scenario 1: $P(\hspace*{0.2em}|$sample log OR*$| > |$observed log OR$| = 1.184\hspace*{0.2em}) \leq
0.006$\\
Scenario 2: $P(\hspace*{0.2em}|$sample log OR*$| > 1.184\hspace*{0.2em}) > 0.006$\\[1.5ex]
where sample log OR* is the sample log OR in a future hypothetical sample of the same size as
the observed sample.
Under the assumptions that have been made, Scenario~1 would be true if and only if log OR lies in
the interval $[-0.1, 0.1]$, while Scenario~2 would be true if and only if log OR does not lie in
this interval.
Given that if a flat improper prior was placed over all values of log OR, the posterior
probability of log OR lying in the interval $[-0.1,0.1]$ would be 0.004 (as given in row~1 of
Table~1), we would expect that the probability that is assigned to Scenario~1 would be
substantially higher than 0.004.

In rows~7 and~8 of Table~1, a high value (0.05) and a low value (0.01) is proposed for the
probability that could be assigned to Scenario~1, which, under the assumptions that have
been made, is equal to the post-data probability that log OR lies in the interval $[-0.1,0.1]$.
The 95\% central post-data intervals for the odds ratio presented in these two rows have been
deduced from these two proposed probability values by applying the same assumptions as were relied
upon in earlier calculations.
It can be seen that these intervals for the odds ratio are very similar to the post-data intervals
for the odds ratio given in rows~4 and~6 that correspond to assigning the same high value (0.05)
and the same low value (0.01) to the probability that log OR $\geq -0.1$ as part of the one-sided P
value method.
\vspace*{2.5ex}

\subsection{Second new study}

Let us again imagine that the new data that needs to be analysed are event rates from the CLARIFY
study, but this time those that relate to the effect of the treatment compared to the control on
the occurrence of heart attacks only. The event counts for this data set can be found in the legend
of Table~2. The sample odds ratio for these event counts is 0.311 and, under the same assumptions
as made earlier, the 95\% confidence interval for the odds ratio is $(0.106, 0.913)$, as indicated
in row~7 of Figure~2.
Clearly this interval does not contain the null effect (odds ratio = 1) but, in contrast to the
previous example, the upper limit of this confidence interval is not that far away from the null
effect.

\begin{table}[tb]
\begin{flushleft}
{\small {\bf Table 2. Analysis of data from the CLARIFY study for myocardial infarction.}}\\[0.5ex]
{\footnotesize Data summary: $n_t = 74$, $n_c = 74$, $e_t = 5$, $e_c = 14$.}
\end{flushleft}
\renewcommand{\arraystretch}{1.1}
\begin{adjustwidth}{-1in}{-1in}
\begin{center}
{\footnotesize
\begin{tabular}{|l|l|l|l|l|l|}
\hline
Method & Prior probability & Significance & Post-data & Post-data & 95\% central\\
& of log OR $\in$ & measure & probability of & probability of & post-data\\
& $[-0.1,0.1]$ & & log \hspace*{-0.1em}OR $\!\in\! [-0.1,0.1]$ & log OR $\geq -0.1$
& interval for OR\\
\hline
Bayesian with flat & \multirow{2}{*}{N/a} & \multirow{2}{*}{N/a} & \multirow{2}{*}{0.0154} &
\multirow{2}{*}{0.026} & \multirow{2}{*}{(0.106, 0.913)}\\
improper prior & & & & & \\
\hline
\multirow{2}{*}{Two-step Bayesian}
& 0.5* & \multirow{2}{*}{N/a} & 0.232 & 0.242 & (0.112, 1.082)\\
& 0.8* & & 0.547 & 0.557 & (0.128, 1.093)\\
\hline
\multirow{3}{28mm}{One-sided P value hybrid} & \multirow{3}{*}{N/a} &
\multirow{3}{23mm}{One-sided P value = 0.0259} & 0.189 & 0.20* & (0.111, 1.078)\\
& & & 0.089 & 0.10* & (0.108, 1.054)\\
& & & 0.039 & 0.05* & (0.106, 1.007)\\
\hline
\multirow{2}{28mm}{Q value hybrid} & \multirow{2}{*}{N/a} & \multirow{2}{23mm}{Q value = 0.0365}
& 0.20* & 0.211 & (0.111, 1.079)\\
& & & 0.05* & 0.061 & (0.107, 1.023)\\
\hline
\end{tabular}}
\end{center}
\end{adjustwidth}
\vspace*{-2ex}
\begin{flushleft}
{\footnotesize The asterisks (*) indicate probabilities determined by the analyst.}
\end{flushleft}
\vspace*{2ex}
\end{table}

The results in Table~2 correspond to applying exactly the same methods and assumptions that were
used to generate the results in Table~1. The one-sided P value of 0.0259 and the Q value of 0.0365
calculated on the basis of the current data set, moderately rather than strongly disfavour the
possibility of log OR lying in the interval $[-0.1,0.1]$. For this reason, the post-data
probabilities that have been proposed in Table~2 for the event that log~OR $\geq -0.1$ in the case
of the one-sided P value method, and the event that log~OR $\in [-0.1,0.1]$ in the case of the Q
value method, are small to moderate probabilities, that is, the values 0.05, 0.10 and 0.20.
It can be seen from the last column of Table~2 that, in contrast to the 95\% confidence interval
for the odds ratio, all solutions proposed in this table for the three methods of primary
interest, i.e.\ the two-step Bayesian method, the one-sided P value method and the Q value method,
correspond to 95\% central post-data intervals for the odds ratio that contain the null effect.
\vspace*{3ex}

\subsection{Third new study}

Finally, let us imagine that the new data that needs to be analysed are event rates from the
STAMINA study, which like the CLARIFY study is a small clinical trial included in the meta-analysis
that appears in Figure~4 of Andraws, Berger and Brown~(2005), with the medical outcome being the
same as in this meta-analysis, i.e.\ acute coronary syndromes.
The events counts for this data set, which can be found in the legend of Table~3, give rise to a
sample odds ratio of 0.640 and, making the same assumptions as previously, the 95\% confidence
interval for the odds ratio is $(0.344, 1.192)$, as indicated in the last row of Figure~2.
In contrast to the previous examples, it can be seen that this confidence interval contains the
null effect.

\begin{table}[tb]
\begin{flushleft}
{\small {\bf Table 3. Analysis of data from the STAMINA study for acute coronary
syndromes.}}\\[0.5ex]
{\footnotesize Data summary: $n_t = 111$, $n_c = 107$, $e_t = 23$, $e_c = 31$.} 
\end{flushleft}
\renewcommand{\arraystretch}{1.1}
\begin{adjustwidth}{-1in}{-1in}
\begin{center}
{\footnotesize
\begin{tabular}{|l|l|l|l|l|l|}
\hline
Method & Prior probability & Significance & Post-data & Post-data & 95\% central\\
& of log OR $\in$ & measure & probability of & probability of & post-data\\
& $[-0.1,0.1]$ & & log \hspace*{-0.1em}OR $\!\in\! [-0.1,0.1]$ & log OR $\geq -0.1$
& interval for OR\\
\hline
Bayesian with flat & \multirow{2}{*}{N/a} & \multirow{2}{*}{N/a} & \multirow{2}{*}{0.095} &
\multirow{2}{*}{0.138} & \multirow{2}{*}{(0.344, 1.192)}\\
improper prior & & & & & \\
\hline
\multirow{2}{*}{Two-step Bayesian}
& 0.5* & \multirow{2}{*}{N/a} & 0.446 & 0.488 & (0.369, 1.112)\\
& 0.8* & & 0.763 & 0.805 & (0.423, 1.099)\\
\hline
One-sided P value & \multirow{3}{*}{N/a} & One-sided P value &
\multirow{3}{*}{\begin{tabular}{@{}l} 0.5*\\ 0.8* \end{tabular}} &
\multirow{3}{*}{\begin{tabular}{@{}l} 0.543\\ 0.843 \end{tabular}} &
\multirow{3}{*}{\begin{tabular}{@{}l}(0.375, 1.104)\\ (0.437, 1.098) \end{tabular}}\\
hybrid / Q value & & = 0.1379 / & & &\\
hybrid & & Q value = 0.1805 & & &\\
\hline
\end{tabular}}
\end{center}
\end{adjustwidth}
\vspace*{-2ex}
\begin{flushleft}
{\footnotesize The asterisks (*) indicate probabilities determined by the analyst.}
\end{flushleft}
\vspace*{2ex}
\end{table}

The results in Table~3 correspond to applying exactly the same methods and assumptions that were
used to generate the results in the earlier tables. The one-sided P value of 0.1379 and the Q value
of 0.1805 calculated on the basis of the data set being considered, at most only marginally
disfavour the possibility of log OR lying in the interval $[-0.1,0.1]$. Therefore, if there was a
substantial belief that log~OR lay in this interval before the data was observed, there still
should be a substantial, albeit slightly lower, belief in this possibility after the data has been
observed.

Since in the two-step Bayesian method, a substantial level of prior belief that log OR lies in the
interval $[-0.1,0.1]$ has been approximately modelled by assigning probabilities of 0.5 and 0.8 to
this event, it has been proposed in Table~3 that these two values approximately represent
probabilities that, under the one-sided P value method or the Q value method, could be reasonably
assigned to this event after the data has been observed. In doing this, we should note that this is
not a strict application of the one-sided P value method as post-data probabilities are being
directly assigned by the analyst to the event $\{\hspace*{0.2em}\mbox{log OR} \in
[-0.1,0.1]\hspace*{0.1em}\}$ rather than the event $\{\hspace*{0.2em}\mbox{log OR} \geq
-0.1\hspace*{0.05em}\}$. Nevertheless, this would seem to be an acceptable modification in this
particular case.

It can be seen from the last column of Table~3 that all solutions proposed in this table for the
methods of primary interest, i.e.\ those reported from the second row onwards, correspond to 95\%
central post-data intervals for the odds ratio that are substantially narrower than the 95\%
confidence interval for the odds ratio and, indeed, are contained within this confidence interval.
\vspace*{3ex}

\section{Discussion}

At the end of Section~\ref{sec3.3}, a comparison was made between the two hybrid approaches
put forward in the present paper, i.e.\ the one-sided P value method and the Q value method. To
complete an assessment of the advantages and disadvantages of all the approaches that have been
proposed to tackle the problem of interest, a two-way comparison will now be made between the
two-step Bayesian approach on the one hand, and these hybrid approaches on the other.

\vspace*{3ex}
\noindent
{\bf Sensitivity to the size of {\large \boldmath $\varepsilon$}} 

\vspace*{1ex}
\noindent
It was set down as a requirement in the Introduction that inferences made about the parameter of
interest $\theta$ should be, to some degree, insensitive to the width of the special interval for
$\theta$, or in other words, to the value of $\varepsilon$. The two-step Bayesian approach fulfils
this requirement if the height of the likelihood function for $\theta$ changes relatively little
over the interval $[\theta_0-\varepsilon, \theta_0+\varepsilon]$ and in the immediate vicinity of
this interval. The hybrid approaches, on the other hand, fulfil the same requirement if the width
of the interval $[\theta_0-\varepsilon, \theta_0+\varepsilon]$ is relatively small compared to the
variance of the sampling distribution for the estimator $\hat{\theta}$. Although, in any given
practical situation it may be difficult to exactly specify the size of $\varepsilon$, we know of
course that $\varepsilon$ must be reasonably small and, in general, this will be sufficient to
ensure that the two-step Bayesian and hybrid approaches will be fairly insensitive to even quite
large changes to the relative size of $\varepsilon$.

\vspace*{3ex}
\noindent
{\bf Bayesian coherency}

\vspace*{1ex}
\noindent
As mentioned in Section~\ref{sec2.2}, the two-step Bayesian approach is incoherent from a Bayesian
perspective due to the prior density for $\theta$ that is used in the first step of the method
being different from the one used in the second step.
The one-sided P value and the Q value
methods, on the other hand, can be regarded as being coherent in terms of the limited way in which
they use Bayesian inference, that is, using the data to update a prior density for $\theta$ 
that is conditioned on $\theta$ lying outside of the interval $[\theta_0-\varepsilon,
\theta_0+\varepsilon]$, or conditioned on $\theta$ lying inside this interval, to a posterior
density that is conditioned in the same way.
However, if the post-data density for $\theta$ over all values of $\theta$ produced by either of
these two hybrid approaches is interpreted from a Bayesian point of view as being a posterior
density for $\theta$, then in general it will not be possible to regard this unconditional density
as being coherent, as only in special cases will it correspond to a given unconditional prior
density for $\theta$.

\vspace*{3ex}
\noindent
{\bf Degree of subjectivity}

\vspace*{1ex}
\noindent
It seems advantageous that using the hybrid approaches, a true post-data probability can be
assigned to the event of $\theta$ lying in the interval $[\theta_0-\varepsilon, \theta_0 +
\varepsilon]$, rather than only being able to determine a lower bound for this probability as is
the case when using the two-step Bayesian approach. However, with all these approaches there is a
large element of subjectivity in how either a true probability, or lower bound probability, for the
event in question is determined. In particular, in the case of the hybrid approaches, the
determination of the probability of this event depends on the interpretation of a one-sided P value
or a Q value through the comparison of the type of opposing scenarios discussed in
Sections~\ref{sec3.2} and~\ref{sec3.3}. By contrast, in the case of the two-step Bayesian approach,
the construction of lower bounds on the probability that $\theta \in [\theta_0-\varepsilon,
\theta_0 + \varepsilon]$ will be sensitive to the choice made for the class of prior densities for
$\theta$ over which such lower bounds are valid, e.g.\ whether this class of prior densities, when
conditioned on $\theta$ not lying in $[\theta_0 - \varepsilon, \theta_0+\varepsilon]$, is more
general or more specific than the class of densities $G_{\mbox{\scriptsize NIS}}$.

\vspace*{3ex}
\noindent
{\bf Ability to corroborate sharp hypotheses}

\vspace*{1ex}
\noindent
An advantage that may be associated with the two-step Bayesian approach is that, in certain
circumstances, observing the data can augment the level of belief that the parameter $\theta$ lies
in the interval $[\theta_0-\varepsilon, \theta_0+\varepsilon]$.
As mentioned in Section~\ref{sec3.2}, this is not naturally facilitated by the one-sided P value
method, which is a comment that also applies to the Q value method.
However, the extent to which, via the use of the two-step Bayesian approach, the data can
corroborate the hypothesis that $\theta$ lies in the interval $[\theta_0-\varepsilon,
\theta_0+\varepsilon]$ may be quite limited.

More specifically, if the likelihood function is a unimodal function with its maximum lying in the
interval $[\theta_0-\varepsilon, \theta_0+\varepsilon]$ and, as in Section~\ref{sec4.2}, the prior
density for $\theta$ conditional on $\theta$ lying in this interval, i.e.\ the density $h(\theta)$,
is assumed to be a uniform density over this interval, then the lower limit on the posterior
probability of $\theta$ lying in the interval $[\theta_0-\varepsilon, \theta_0+\varepsilon]$ over
all $g(x) \in G_{\mbox{\scriptsize NIS}}$ will be
\begin{equation}
\label{equ7}
p_{L}(\theta \in [\theta_0 - \varepsilon, \theta_0 + \varepsilon]\,|\,x) = 
\alpha (M_{\mbox{\footnotesize \tt inside}} / M_{\mbox{\footnotesize \tt limits}})
\end{equation}
where $M_{\mbox{\footnotesize \tt inside}}$ is the mean height of the likelihood function over the
interval $[\theta_0-\varepsilon, \theta_0+\varepsilon]$, $M_{\mbox{\footnotesize \tt limits}}$ is
the mean of the height of the likelihood function at the two limits of this interval and, as a
reminder, $\alpha$ is the prior probability of $\theta$ lying in the interval.
Since it would not be expected, given the nature of the problem of interest, that the height of
the likelihood function will be much lower at the limits of the interval $[\theta_0-\varepsilon,
\theta_0+\varepsilon]$ than at its maximum, which must lie in this interval, the factor
$M_{\mbox{\footnotesize \tt inside}} / M_{\mbox{\footnotesize \tt limits}}$ in
equation~(\ref{equ7}) will not be expected to be much greater than one.
Therefore, in practice, the lower limit on the posterior probability of $\theta$ lying in the
interval $[\theta_0-\varepsilon, \theta_0+\varepsilon]$ will generally not be much greater than the
prior probability $\alpha$ that is placed on this event. On the other hand, this lower limit, in
general, will be less than $\alpha$ if the maximum of the likelihood function does not lie in the
interval $[\theta_0-\varepsilon, \theta_0+\varepsilon]$.

\vspace*{3ex}
\noindent
{\bf Incorporation of nuisance parameters}

\vspace*{1ex}
\noindent
Another advantage that may be associated with the two-step Bayesian approach is that, in general,
it may be easier to extend this approach to deal with the presence of nuisance parameters
than is the case with the one-sided P value or the Q value approaches. Using the standard Bayesian
paradigm, the most logical way of incorporating the presence of nuisance parameters into the
two-step Bayesian approach outlined in Section~\ref{sec2.2} is to specify a (joint) prior density
for the nuisance parameter(s) conditional on the parameter of interest $\theta$. The calculations
of posterior probabilities and densities in the first and second steps of this approach are then
carried out using this conditional prior density and a class of marginal prior densities for the
parameter $\theta$, e.g.\ the class of densities $G_{\mbox{\scriptsize NIS}}$.
However, there may be practical difficulties in applying such a multivariate approach.
In particular, it may be difficult to find a joint density for the nuisance parameters conditional
on the parameter $\theta$ that adequately represents prior opinion about the dependency between all
these parameters, and also there is the potentially very challenging task of finding the lower
limit for the posterior probability of $\theta$ lying in the interval $[\theta_0-\varepsilon,
\theta_0+\varepsilon]$ over a class of joint prior densities for the parameter $\theta$ and the
nuisance parameters.

For the one-sided P value and the Q value methods, extensions to cases with nuisance parameters
could follow along similar lines to work that has already been done in extending hypothesis testing
methods to deal with the presence of such parameters. In particular, the hypothesis testing methods
that have been developed to tackle this issue could be categorized into conventional methods, e.g.\
the work of Tsui and Weerahandi~(1989) and Berger and Boos~(1994), bootstrap methods, e.g.\ Davison
and Hinkley~(1997), and Bayesian methods, e.g.\ Box~(1980) and Meng~(1994). However, it could be
argued that this is a controversial area of research, as is indicated by the large number of
competing methods that have been proposed to resolve the issue concerned, and the lack of
consensus in the choice amongst these methods.

\vspace*{3ex}
\noindent
{\bf Representation of prior beliefs}

\vspace*{1ex}
\noindent
To finish this discussion, we will address the issue of whether the lower bounds on the posterior
probability of the event of the parameter $\theta$ lying in the interval $[\theta_0-\varepsilon,
\theta_0+\varepsilon]$ obtained in the second step of the two-step Bayesian approach can invalidate
solutions to the problem of interest obtained using the hybrid approaches of Section~\ref{sec3}.
In particular, it seems relevant to inquire as to whether post-data probabilities for this event
obtained using the one-sided P value, or the Q value methods, can be valid if they are smaller than
minimum bounds on these probabilities established by performing a Bayesian sensitivity analysis
over a very general class of prior densities for $\theta$.

The answer to this question can be found at the heart of a long-standing debate concerning the
foundations of statistical inference. In this regard, we should begin by taking into account that,
in the Bayesian approach to inference, it is assumed that beliefs about a parameter, before the
data was observed, can be adequately represented by placing a probability distribution over this
parameter, i.e.\ a prior distribution.
This assumption clearly remains even if a Bayesian sensitivity analysis is performed over a class
of such prior distributions. However in the context of interest, first, it should not be taken for
granted that a substantial degree of prior belief that $\theta$ lies in the narrow interval
$[\theta_0 -  \varepsilon, \theta_0+\varepsilon]$ can be adequately represented by assigning a
probability value to this event.
Second, it is even more debatable that prior opinion about $\theta$ conditional on this event not
occurring can be represented satisfactorily by placing a (proper) probability distribution over all
values for $\theta$ not lying in the interval $[\theta_0 - \varepsilon, \theta_0 + \varepsilon]$.

The need to express prior beliefs using probabilities and probability distributions in the
Bayesian approach places a rigid structure on the way such beliefs can be represented, which may
lead to prior opinion being excessively favoured over information contained in the data.
On the other hand, it could be argued that the one-sided P value and the Q value methods allow
prior beliefs to be incorporated into the inferential process in a more flexible way, with the
result that these beliefs can be modified to a greater extent if they appear to be inconsistent
with the data. In conclusion, it is evident that the two-step Bayesian approach does not
necessarily stand in contradiction to the hybrid approaches of Section~3 from a foundational
viewpoint.
Therefore, the choice between the approaches that have been put forward should be made primarily on
the basis of all the other substantive issues that have been highlighted throughout the present
paper.

\pagebreak

\noindent
{\bf References}
\vspace*{1ex}

\begin{description}

\item[] Andraws, R., Berger, J. S. and Brown, D. L. (2005). Effects of antibiotic therapy on
outcomes of patients with coronary artery disease: a meta-analysis of randomized controlled trials.
\emph{Journal of the American Medical Association}, {\bf 293}, 2641--2647.

\item[] Anscombe, F. J. (1990). The summarizing of clinical experiments by significance levels.
\emph{Statistics in Medicine}, {\bf 9}, 703--708.

\item[] Berger, J. O. and Delampady, M. (1987). Testing precise hypotheses. \emph{Statistical
Science}, {\bf 2}, 317--335.

\item[] Berger, J. O. and Sellke, T. (1987). Testing a point null hypothesis: the irreconcilability
of P values and evidence. \emph{Journal of the American Statistical Association}, {\bf 82},
112--122.

\item[] Berger, R. L. and Boos, D. D. (1994). P values maximized over a confidence set for the
nuisance parameter, \emph{Journal of the American Statistical Association}, {\bf 89}, 1012--1016.

\item[] Box, G. E. P. (1980). Sampling and Bayes' inference in scientific modelling and robustness
(with discussion). \emph{Journal of the Royal Statistical Society, Series A}, {\bf 143}, 383--430.

\item[] Davison, A. C. and Hinkley D. V. (1997). \emph{Bootstrap methods and their application},
Cambridge University Press, Cambridge.

\item[] Delampady, M. (1989). Lower bounds on Bayes factors for interval null hypotheses.
\emph{Journal of the American Statistical Association}, {\bf 84}, 120--124.

\item[] DerSimonian, R. and Laird, N. (1986). Meta-analysis in clinical trials. \emph{Controlled
Clinical Trials}, {\bf 7}, 177--188.

\item[] Edwards, W. (1965). Tactical note on the relation between scientific and statistical
hypotheses. \emph{Psychological Bulletin}, {\bf 63}, 400--402.

\item[] Edwards, W., Lindman, H. and Savage, L. J. (1963). Bayesian statistical inference for
psychological research. \emph{Psychological Review}, {\bf 70}, 193--242.

\item[] Greenland, S. (2011). Null misinterpretation in statistical testing and its impact on
health risk assessment. \emph{Preventive Medicine}, {\bf 53}, 225--228.

\item[] Meng, X. L. (1994). Posterior predictive p-values. \emph{Annals of Statistics},
{\bf 22}, 1142--1160.

\item[] Oakes, M. W. (1990). \emph{Statistical Inference}. Epidemiology Resources Inc., Chestnut
Hill.

\item[] Royall, R. M. (1986). The effect of sample size on the meaning of significance tests.
\emph{The American Statistician}, {\bf 40}, 313--315.

\item[] Tsui, K. W. and Weerahandi, S. (1989). Generalized p-values in significance testing of
hypotheses in the presence of nuisance parameters. \emph{Journal of the American Statistical
Association}, {\bf 84}, 602--607.

\item[] Woolf, B. (1955). On estimating the relationship between blood group and disease.
\emph{Annals of Human Genetics}, {\bf 19}, 251--253.

\end{description}

\end{document}